\documentclass[10pt]{article}
\hoffset - 20 mm
\usepackage{amsmath,amssymb,amsthm}
\usepackage{graphicx}
\usepackage{color}
\catcode`\@=11 \@addtoreset{equation}{section}

\catcode`\@=12

\newtheorem{Theorem}{Theorem}[section]
\newtheorem{Proposition}{Proposition}[section]
\newtheorem{Lemma}{Lemma}[section]
\newtheorem{Corollary}{Corollary}[section]
\newtheorem{Remark}{Remark}[section]

\newcommand{\bTheorem}[1]{
\begin{Theorem} \label{T#1} }
\newcommand{\eT}{\end{Theorem}}

\newcommand{\bProposition}[1]{
\begin{Proposition} \label{P#1}}
\newcommand{\eP}{\end{Proposition}}

\newcommand{\bLemma}[1]{
\begin{Lemma} \label{L#1} }
\newcommand{\eL}{\end{Lemma}}

\newcommand{\bCorollary}[1]{
\begin{Corollary} \label{C#1} }
\newcommand{\eC}{\end{Corollary}}

\newcommand{\bFormula}[1]{
\begin{equation} \label{#1}}
\newcommand{\eF}{\end{equation}}

\newcommand{\vr}{\varrho}

\newcommand{\vt}{\vartheta}

\newcommand{\vu}{\vc{u}}
\newcommand{\vc}[1]{{\vec #1}}
\newcommand{\Div}{{\rm div}_x}
\newcommand{\Grad}{\nabla_x}
\newcommand{\Curl}{{\rm curl}_x}

\font\F=msbm10 scaled 1000
\newcommand{\R}{\mbox{\F R}}

\date{}
\makeindex
\begin{document}

\title{Global existence of a radiative Euler system coupled to an electromagnetic field}
\author{X. Blanc, B. Ducomet, \v S. Ne\v casov\' a}
\maketitle

 \begin{abstract}
 \noindent We study the Cauchy problem for a system of equations corresponding to a singular limit of radiative hydrodynamics,
 namely the 3D radiative compressible Euler system coupled to an electromagnetic field.
Assuming smallness hypotheses for the data, we prove that the problem admits a unique global smooth solution and study its
asymptotics.
\end{abstract}



\section{Introduction}

In \cite{BDN}, after the studies of Lowrie, Morel and Hittinger \cite{LMH} and Buet and Despr\'es  \cite{BD} we considered a
 singular limit for a compressible inviscid radiative flow where the motion of the fluid is given by the Euler system for the
evolution of the density $\vr = \vr (t,x)$, the velocity
field $\vu = \vu(t,x)$, and the absolute temperature $\vt =\vt(t,x)$, and where radiation is described in the limit by
 an extra temperature $T_r=T_r(t,x)$. All of these quantities are functions of the time $t$ and the Eulerian spatial
coordinate $x \in \R^3$.

 In \cite{BDN} we proved that the associated Cauchy problem admits a unique global smooth solution, provided that the data
are small enough perturbations of a constant state.

 In \cite{BDN2} we coupled the previous model to the electromagnetic field
through the so called magnetohydrodynamic (MHD) approximation, in presence of thermal and radiative dissipation.
 Hereafter we consider the perfect non-isentropic Euler-Maxwell's system and we also consider a radiative coupling through
a pure convective transport equation for the radiation (without dissipation). Then we deal with a pure hyperbolic system
with partial relaxation (damping on velocity).

More specifically the system of equations to be studied for the unknowns $(\vr,\vu,\vt,E_r,\vec B, \vec E)$ reads
\bFormula{j1bis}
\partial_t \vr + \Div (\vr \vu) = 0,
\eF
\bFormula{j2bis}
\partial_t (\vr \vu) + \Div (\vr \vu \otimes \vu) + \Grad (p+p_r) = -\rho\left(\vec E+\vc u\times\vec B\right)-\nu\rho\vu,
\eF
\bFormula{j3bis}
\partial_t \left( \vr E \right)
+ \Div \left( (\vr E+p)  \vu\right) +  \vu \cdot\Grad p_r= - \sigma_a\left(a\vt^4-E_r\right)-\rho\vec E\cdot\vc u,
\eF
\bFormula{j4bis}
\partial_t E_r
+ \Div \left( E_r  \vu\right) + p_r\Div \vu  = -\sigma_a\left(E_r-a\vt^4\right),
\eF
\bFormula{j5bis}
\partial_t \vc{B} + \Curl \vc{E} = 0,\
\eF
\bFormula{j6bis}
\partial_t \vc{E} - \Curl \vc{B}  = \vr\vec u,
\eF
\bFormula{j7bis}
\Div \vc{B} = 0,\
\eF
\bFormula{j8bis}
\Div \vc{E}  = \overline\vr-\vr,
\eF
where $\vr$ is the density, $\vu$ the velocity, $\vt$ the temperature of matter, $E= \frac{1}{2} |\vu|^2 + e(\vr,\vt)$ is the total mechanical energy,
 $E_r$ is the radiative energy related to the temperature of radiation $T_r$ by $E_r=aT_r^4$ and $p_r$ is the radiative pressure given
 by $p_r=\frac{1}{3} aT_r^4=\frac{1}{3}\ E_r$, with $a>0$. Finally $\vec E$ is the electric field and
 $\vc{B}$ is the magnetic induction,

 We assume that the pressure $p(\vr,\vt)$ and the internal energy $e(\vr,\vt)$ are positive smooth functions of their arguments with
\[
C_v:=\frac{\partial e}{\partial\vt}>0,\ \ \frac{\partial p}{\partial\vr}>0,
\]
 and we also suppose for simplicity that $\nu=\frac{1}{\tau}$ (where $\tau>0$ is a momentum-relaxation time), $\mu,\sigma_a$ and $a$ are positive constants.

A simplification appears if one observes that, provided that equations (\ref{j7bis}) and (\ref{j8bis}) are satisfied at $t=0$,
they are satisfied for any time $t>0$ and consequently they can be discarded from the analysis below.

Notice that the reduced system (\ref{j1bis})-(\ref{j4bis}) is the non equilibrium regime of radiation hydrodynamics
 introduced by Lowrie, Morel and Hittinger \cite {LMH} and more recently by Buet and Despr\' es \cite{BD}, and studied mathematically
 by Blanc, Ducomet and  Ne\v casov\' a \cite{BDN}. Extending this last analysis, our goal in this work is to prove
global existence of solutions for the system (\ref{j1bis}) - (\ref{j8bis})
 when data are sufficiently close to an equilibrium state, and study their large time behaviour.

Just mention for completeness that related non isentropic Euler-Maxwell systems have been the object of a number of studies in the recent past.
 Let us quote some recent works:
Y. Feng, S. Wang, S. Kawashima \cite{FWK},
Y. Feng, S. Wang, X. Li \cite{FWL},
J.W. Jerome \cite{J},
C. Lin, T. Goudon \cite{LG},
Z. Tan, Y. Wang \cite{TW}
and
J. Xiu, J. Xiong \cite{XX}.

\vskip0.25cm
In the following we show that the ideas used by Y. Ueda, S. Wang and S. Kawashima in \cite{UWK} \cite{UK} in the isentropic case can be extended
 to the (radiative) non isentropic system (\ref{j1bis}-\ref{j6bis}). To this purpose we follow the following plan: in Section 2 we present the main results, then (Section 3) we prove well-posedness of system (\ref{j1bis}-\ref{j6bis}).
Finally in Section 4 we prove the large time asymptotics of the solution.

\section{Main results}

We are going to prove that system \eqref{j1bis}-\eqref{j8bis} has a global smooth solution close to any equilibrium
state. Namely we have

\begin{Theorem}
\label{th:existence_Maxwell}
  Let $\left(\overline \vr, 0, \overline\vt,\overline{E_r},\overline{\vec B},0\right)$ be a constant
  state with $\overline \vr >0$, $\overline\vt >0$ and $\overline{E_r}>0$ with compatibility condition $\overline{E_r}=a\overline\vt^4$ and suppose that $d\geq3$.

There exists $\varepsilon>0$ such that, for any initial state $\left(\vr_0,\vu_0,\vt_0,E_r^0,\vec B_0,\vec E_0\right)$ satisfying
\[
 \Div \vec B_0= \vr_0-\overline\vr,\ \ \ \Div \vec B_0=0,
\]
\[
\left(\vr_0-\overline \vr, \vu_0, \vt_0-\overline\vt,{E_r}_0-\overline{E_r},\vec B_0-\overline{\vec B},\vec E_0\right)\in H^d,
\]
and
\begin{equation}
  \label{eq:donnees_petites_maxwell}
  \left\|\left(\vr_0,\vu_0,\vt_0,E_r^0,\vec B_0,\vec E_0\right) - \left(\overline \vr, 0, \overline\vt,\overline{E_r},\overline{\vec B},0\right) \right\|_{H^d}
\leq \varepsilon,
\end{equation}
there exists a unique global solution $\left(\vr,\vu,\vt,E_r,\vec B,\vec E\right)$ to
\eqref{j1bis}-\eqref{j8bis}, such that
\[
\left(\vr-\overline\vr,\vu,\vt-\overline \vt,E_r-\overline{E_r},\vec B-\overline{\vec B},\vec E\right) \in
C\left([0,+\infty);H^d \right)\cap C^1\left([0,+\infty);H^{d-1} \right).
\]
 In addition, this solution satisfies the following energy inequality:
\[
  \left\|
(\vr-\overline\vr,\vu,\vt-\overline \vt,E_r-\overline{E_r},\vec B-\overline{\vec B},\vec E)(t)
\right\|_{H^d}
\]
\[
 + \int_0^t
\left(
 \left\|
\left(\vr-\overline\vr,\vu,\vt-\overline \vt,E_r-\overline{E_r}\right)(\tau)
\right\|^2_{H^{d}}
 + \left\|\Grad \vec B(\tau)\right\|^2_{H^{d-2}}
 + \left\|\vec E(\tau)\right\|^2_{H^{d-1}}
\right)\ d\tau
\]
\begin{equation}
\leq C \left\|\left(\vr_0-\overline\vr,0,\vt_0-\overline \vt,E_r^0-\overline{E_r},\vec B_0-\overline{\vec B},\vec E_0\right)\right\|_{H^d}^2,
\label{eq:energy_estimate_maxwell}
\end{equation}
for some constant $C>0$ which does not depend on $t$.
\end{Theorem}
The large time behaviour of the solution is described as follows

\begin{Theorem}
\label{th:asymptotics_Maxwell}
Let $d\geq 3$.

The unique global solution $\left(\vr,\vu,\vt,E_r,\vec B,\vec E\right)$ to
\eqref{j1bis}-\eqref{j8bis} defined in Theorem \ref{th:existence_Maxwell} converges to the constant state
$\left(\overline\vr,\vec 0,\overline \vt,\overline{E_r},\overline{\vec B},\vec 0\right)$ uniformly in $x\in{\mathbb R}^3$ as $t\to\infty$.
More precisely
\begin{equation}
\left\|(\vr-\overline\vr,\vu,\vt-\overline \vt,E_r-\overline{E_r},\vec E)(t)\right\|_{W^{d-2,\infty}}\to 0\ \ \ \mbox{as}\ t\to\infty.
\label{as1}
\end{equation}
Moreover if $d\geq 4$
\begin{equation}
\left\|(\vec B-\overline{\vec B})(t)\right\|_{W^{d-4,\infty}}\to 0\ \ \ \mbox{as}\ t\to\infty,
\label{as2}
\end{equation}
\end{Theorem}

\begin{Remark}
\label{r1}
Note that, due to lack of
 dissipation by viscous, thermal and radiative fluxes,
 the Kawashima-Shizuta stability criterion (see \cite{SK} and \cite{BZ}) is not satisfied for the system under study and techniques of \cite{K} relying on the existence of
a compensating matrix do not apply. However we will check that radiative sources play the role of relaxation terms for temperature and radiative energy and
 will lead to global existence for the system.
\end{Remark}

\section{Global existence}
\label{sec:linearized-EM}

\subsection{A priori estimates}

Multiplying (\ref{j2bis}) by $\vu$, (\ref{j5bis}) by $\vec B$, (\ref{j6bis}) by $\vec E$ and adding the result
to equations (\ref{j3bis}) and (\ref{j4bis}) we get the total energy conservation law
\bFormula{ienerEM}
\partial_t \left( \frac{1}{2}\vr \left|\vu\right|^2+\vr e +E_r+\frac{1}{2}\left( |\vec B|^2+ |\vec E|^2\right)\right)
+ \Div \left( (\vr E+E_r)  \vu + (p+ p_r)\vu+ \vec E\times\vec B\right)=0.
\eF
Introducing the entropy $s$ of the fluid by the Gibbs law $\vt ds=de+pd\left(\frac{1}{\vr}\right)$ and denoting by $S_r:=\frac{4}{3}aT_r^3$ the radiative entropy,
 equation (\ref{j4bis}) rewrites
\bFormula{iEREM}
\partial_t S_r
+ \Div \left( S_r  \vu\right)  =
-\sigma_a\frac{E_r-a\vt^4}{T_r}.
\eF
The internal energy equation is
\begin{equation}
\partial_t (\vr e) + \Div (\vr e \vu) + p\Div \vu - \nu\vr \left|\vu\right|^2
=
- \sigma_a\left(a\vt^4-E_r\right),
\label{internEM}
\end{equation}
and dividing it by $\vt$, we get the entropy equation for matter
\begin{equation}
\partial_t (\vr s) + \Div (\vr s \vu) -\frac\nu\vt\left|\vu\right|^2 =
 - \sigma_a\frac{a\vt^4-E_r}{\vt}.
\label{entroEM}
\end{equation}
So adding (\ref{entroEM}) and (\ref{iEREM}) we obtain
\bFormula{ientEM}
\partial_t \left( \vr  s +S_r\right)+ \Div \left( (\vr  s  +S_r)\vu\right)
=
\frac{a\sigma_a}{\vt T_r}\left(\vt-T_r\right)^2\left(\vt+T_r\right)\left(\vt^2+T_r^2\right)+ \frac\nu\vt \left|\vu\right|^2.
\eF
Subtracting (\ref{ientEM}) from (\ref{ienerEM}) and using the conservation of mass, we get
\[
\partial_t \left( \frac{1}{2}\vr \left|\vu\right|^2
+H_{\overline{\vt}}(\vr,\vt)-(\vr-\overline{\vr})\partial_{\vr}H_{\overline{\vt}}(\overline{\vr},\overline{\vt})-H_{\overline{\vt}}(\overline{\vr},\overline{\vt})
+H_{r,\overline{\vt}}(T_r)+\frac{1}{2}\left( |\vec B-\overline{\vec B}|^2+ |\vec E|^2\right)
 \right)
\]
\bFormula{Ienerbis}
= \Div \left( (\vr E+E_r)  \vu + (p+ p_r)\vu +\overline{\vt} (\vr s+S_r)  \vu\right)
- \overline{\vt}\frac{a\sigma_a}{\vt T_r}\left(\vt-T_r\right)^2\left(\vt+T_r\right)\left(\vt^2+T_r^2\right)- \frac\nu\vt \left|\vu\right|^2.
\eF
Introducing the Helmholtz functions $H_{\overline{\vt}}(\vr,\vt):=\vr\left(e-\overline{\vt}s\right)$ and $H_{r,\overline{\vt}}(T_r):=E_r-\overline{\vt}S_r$, we check that
the quantities $H_{\overline{\vt}}(\vr,\vt)-(\vr-\overline{\vr})\partial_{\vr}H_{\overline{\vt}}(\overline{\vr},\overline{\vt})-H_{\overline{\vt}}(\overline{\vr},\overline{\vt})$ and
$H_{r,\overline{\vt}}(T_r)-H_{r,\overline{\vt}}(\overline{T}_r)$ are non-negative and strictly coercive functions reaching zero minima at the equilibrium state $(\overline{\vr},\overline{\vt},\overline{E}_r)$.

\newtheorem{l1}{Lemma}
\begin{l1}
Let $\overline{\vr}$ and $\overline{\vt} = \overline{T}_r$ be given positive constants.
Let ${\mathcal O}_1$ and ${\mathcal O}_2$ be the sets defined by

\bFormula{es1}
{\mathcal O}_1:=\left\{
(\vr,\vt)\in\R^2\ :\ \frac{\overline{\vr}}{2}<\vr<2\overline{\vr},\ \frac{\overline{\vt}}{2}<\vt<2\overline{\vt},
\right\}.
\eF
\bFormula{es2}
{\mathcal O}_2:=\left\{
 T_r\in\R\ :\ \frac{\overline{T}_r}{2}<T_r<2\overline{T}_r,
\right\}.
\eF
There exist positive constants $C_{1,2}(\overline{\vr},\overline{\vt})$ and $C_{3,4}(\overline{T}_r)$ such that
\begin{enumerate}
\item
\[
C_1\left(|\vr-\overline{\vr}|^2+|\vt-\overline{\vt}|^2\right)
\leq
H_{\overline{\vt}}(\vr,\vt)
-(\vr-\overline{\vr})\partial_{\vr}H_{\overline{\vt}}(\overline{\vr},\overline{\vt})
-H_{\overline{\vt}}(\overline{\vr},\overline{\vt})
\]
\bFormula{ess1}
\leq
C_2\left(|\vr-\overline{\vr}|^2+|\vt-\overline{\vt}|^2\right),
\eF
for all $(\vr,\vt)\in{\mathcal O}_1$,
\item
\bFormula{ess2}
C_3 |T_r-\overline{T}_r|^2
\leq
 H_{r,\overline{\vt}}(T_r)-H_{r,\overline{\vt}}(\overline{T}_r)
\leq
C_4|T_r-\overline{T}_r|^2,
\eF
for all $T_r\in{\mathcal O}_2$.
\end{enumerate}
\label{l1}
\end{l1}

\noindent {\bf Proof:}
\begin{enumerate}
\item Point $1$ is proved in \cite{FEINOV} and we only sketch the proof for convenience. According to the decomposition
$$\vr\to
H_{\overline{\vt}}(\vr,\vt)
-(\vr-\overline{\vr})\partial_{\vr}H_{\overline{\vt}}(\overline{\vr},\overline{\vt})
-H_{\overline{\vt}}(\overline{\vr},\overline{\vt})
={\mathcal F}(\vr)+{\mathcal G}(\vr),$$
where
${\mathcal F}(\vr)=H_{\overline{\vt}}(\vr,\overline{\vt})
-(\vr-\overline{\vr})\partial_{\vr}H_{\overline{\vt}}(\overline{\vr},\overline{\vt})
-H_{\overline{\vt}}(\overline{\vr},\overline{\vt})$ and
 $ {\mathcal G}(\vr)=H_{\overline{\vt}}(\vr,\vt)-H_{\overline{\vt}}(\vr,\overline{\vt})$,
one checks that ${\mathcal F}$ is strictly convex and reaches a zero minimum at $\overline{\vr}$, while $ {\mathcal G}$ is strictly
decreasing for $\vt<\overline{\vt}$ and strictly
increasing for $\vt>\overline{\vt}$, according to the standard thermodynamic stability properties \cite{FEINOV}. Computing
the derivatives of $H_{\overline{\vt}}$ leads directly to the estimate (\ref{ess1}).

\item Point $2$ follows after properties of
$x\to H_{r,\overline{\vt}}(x)-H_{r,\overline{\vt}}(T_r)=ax^3(x-\frac{4}{3}\overline{\vt})+\frac{a}{3}\overline{\vt}^4$. \hfill$\square$
\end{enumerate}
\vskip0.25cm
Using the previous entropy properties, we have the energy estimate

\begin{Proposition}
  \label{pr:linftyh_estimat}
Let the assumptions of Theorem~\ref{th:existence_Maxwell} be
satisfied with $V = \left(\rho,\vu,\vt,E_r,\vec B,\vec E \right)$,
$\overline V = \left(\overline \rho, 0, \overline\vt,\overline{E_r},\overline{\vec B},0  \right)$.
 Consider a solution $(\vr,\vu,\vt,E_r,\vec B,\vec E)$ of system
\eqref{j1bis}-\eqref{j2bis}-\eqref{j3bis} on $[0,t]$, for some $t>0$.
Then, one gets for a constant $C_0>0$
\begin{equation}
  \label{eq:linftyhd_estimat}
 \left\|V(t) -\overline V \right\|_{L^2}^2 + \int_0^t
\left\|\ \vu(\tau)\right\|_{L^2}^2 d\tau
 \leq C_0\left\|V_0 -\overline V \right\|_{L^2}^2.
\end{equation}
\end{Proposition}
\noindent {\bf Proof:} Defining
\begin{equation}
  \label{eq:def_etabis}
  \eta(t,x) = H_{\overline\vt}(\vr,\vt) - \left(\vr - \overline\vr\right)
    \partial_\vr H_{\overline \vt}\left(\overline\vr,\overline\vt\right)
    -H_{\overline \vt}\left(\overline\vr,\overline\vt\right) +
    H_{r,\overline\vt} \left(T_r \right),
\end{equation}
we multiply \eqref{ientEM} by $\overline\vt$, and subtract the result to
\eqref{ienerEM}. Integrating over $[0,t]\times \R^3$, we find
\[
  \int_{\R^3}\frac{1}{2} \vr \left|\vu\right|^2 +\eta(t,x)+\frac{1}{2} |\vec B-\overline{\vec B}|^2+\frac{1}{2} |\vec E|^2 dx
 + \int_0^t \int_{\R^3}
  \frac{\overline{\vt}}{\vt}\nu\left|\vu\right|^2
\]
\[
     \leq \int_{\R^3}\frac12 \vr_0 \left|\vu_0\right|^2(t) +\eta(0,x)+\frac{1}{2} |\vec B_0-\overline{\vec B}|^2+\frac{1}{2} |\vec E_0|^2\ dx.
\]
Applying Lemma~\ref{l1}, we find \eqref{eq:linftyhd_estimat}.
\hfill $\square$

\vskip0.5cm
Defining for any $d\geq 3$ the auxiliary quantities
\[
{\displaystyle
E(t):=\sup_{0\leq\tau\leq t}
 \|(\vr-\overline \vr,\vu,\vec B-\overline{\vec B},\vec E)(\tau)\|_{W^{1,\infty}},
}
\]
\[
{\displaystyle
F(t):=\sup_{0\leq\tau\leq t}\|(V-\overline V)(\tau)\|_{H^d},
}
\]
\[
{\displaystyle
I^2(t):=\int_0^t\left\|(\vr-\overline\vr,\vu,\vt-\overline\vt,E_r-\overline E_r)\right\|_{L^\infty}^2 d\tau,
}
\]
and
\[
{\displaystyle
D^2(t):=\int_0^t\left(\left\|(\vr-\overline\vr,\vu,\vt-\overline\vt,E_r-\overline E_r)(\tau)\right\|_{H^d}^2
+\left\|\vec E(\tau)  \right\|_{H^{d-1}}^2+\left\| \partial_x\vec B(\tau)  \right\|_{H^{d-2}}^2\right) d\tau,
}
\]
 we can bound the spatial derivatives as follows

\begin{Proposition}
  \label{pr:linftyh-estimates_3}
 Assume that the hypotheses of Theorem~\ref{th:existence_Maxwell} are satisfied. Then, we have for a $C_0>0$
\begin{equation}
\label{eq:linftyhd_estimat_3}
 \left\|\partial_x V(t) \right\|_{H^{d-1}}^2
 + \int_0^t \|\partial_x \vu(\tau)\|^2_{H^{d-1}}\ d\tau
\leq C_0 \left\|\partial_x V_0 \right\|_{H^{d-1}}^2
+C_0 \left(E(t)D(t)^2+F(t)I(t)D(t)\right).
  \end{equation}
\end{Proposition}

\noindent {\bf Proof:} Rewriting the system (\ref{j1bis})-(\ref{j6bis}) in the form
\begin{equation}
\left\{
\begin{array}{c}
\partial_t \vr + \vu\cdot\Grad\vr+\vr\Div  \vu = 0,\\\\
\partial_t  \vu + (\vu\cdot\Grad) \vu + \frac{p_\vr}{\vr}\Grad \vr + \frac{p_\vt}{\vr}\Grad \vt+\frac{1}{3a\vr}\Grad E_r
+\vec E+\vu\times\overline{\vec B}+\nu\vu= -\vc u\times\left(\vec B-\overline{\vec B}\right),\\\\
\partial_t \vt+ (\vu\cdot\Grad) \vt +  \frac{\vt p_\vt}{\vr C_v}\Div  \vu= - \frac{\sigma_a}{\vr C_v}\left(a\vt^4-E_r\right),\\\\
\partial_t E_r+ (\vu\cdot\Grad) E_r   + \frac{4}{3}E_r\Div \vu  = -\sigma_a\left(E_r-a\vt^4\right),\\\\
\partial_t \vc{B} + \Curl \vc{E} = 0,\\\\
\partial_t \vc{E} - \Curl \vc{B} -\overline\vr\vec u = (\vr-\overline\vr)\vec u,
\end{array}
\label{sys}
\right.
\end{equation}
and applying $\partial_x^\ell$ to this system, we get

\[
\partial_t(\partial_x^\ell \vr) + (\vu\cdot\Grad)\partial_x^\ell\vr+\vr\Div  \partial_x^\ell\vu = F_1^\ell,
\]
\[
\partial_t(\partial_x^\ell\vu) + (\vu\cdot\Grad) \partial_x^\ell\vu + \frac{p_\vr}{\vr}\Grad \partial_x^\ell\vr + \frac{p_\vt}{\vr}\Grad \partial_x^\ell\vt+\frac{1}{3a\vr}\Grad \partial_x^\ell E_r
+\partial_x^\ell\vec E+\partial_x^\ell\vu\times\overline{\vec B}+\nu\partial_x^\ell\vu= -\partial_x^\ell\left[\vc u\times\left(\vec B-\overline{\vec B}\right)\right]+F_2^\ell,
\]
\[
\partial_t (\partial_x^\ell\vt)+ (\vu\cdot\Grad) \partial_x^\ell\vt +  \frac{\vt p_\vt}{\vr C_v}\Div  \partial_x^\ell\vu
=
 - \partial_x^\ell\left[\frac{\sigma_a}{\vr C_v}\left(a\vt^4-E_r\right)\right]+F_3^\ell,
\]
\[
\partial_t (\partial_x^\ell E_r)+ (\vu\cdot\Grad) \partial_x^\ell E_r   + \frac{4}{3}E_r\Div \partial_x^\ell\vu
=-\partial_x^\ell\left[\sigma_a\left(E_r-a\vt^4\right)\right]+F_4^\ell,
\]
\[
\partial_t (\partial_x^\ell\vc{B}) + \Curl \partial_x^\ell\vc{E} = 0,
\]
\[
\partial_t (\partial_x^\ell\vc{E}) - \Curl \partial_x^\ell\vc{B} -\overline\vr\partial_x^\ell\vec u = \partial_x^\ell[(\vr-\overline\vr)\vec u],
\]
where
\[
F_1^\ell:=-\left[\partial_x^\ell,\vu\cdot\Grad\right]\vu-\left[\partial_x^\ell,\vr\Div\right]\vu,
\]
\[
F_2^\ell:=-\left[\partial_x^\ell,\vu\cdot\Grad\right]\vu-\left[\partial_x^\ell,\frac{p_\vr}{\vr}\Grad\right]\vr
-\left[\partial_x^\ell,\frac{p_\vt}{\vr}\Grad\right]\vt-\left[\partial_x^\ell,\frac{1}{3a\vr}\Grad \right]E_r,
\]
\[
F_3^\ell:=-\left[\partial_x^\ell,\vu\cdot\Grad\right]\vt-\left[\partial_x^\ell,\frac{\vt p_\vt}{\vr C_v}\Div\right]\vu,
\]
and
\[
F_4^\ell:=-\left[\partial_x^\ell,\vu\cdot\Grad\right]E_r-\left[\partial_x^\ell,\frac{4}{3}E_r\Div\right]\vu.
\]

Then taking the scalar product of each of the previous equations respectively by
$\frac{p_\vr}{\vr^2}\partial_x^\ell \vr,\
\partial_x^\ell \vu,\
\frac{C_v}{\vt}\partial_x^\ell \vt,\
\frac{1}{4a\vr E_r}\partial_x^\ell E_r,\
\partial_x^\ell\vc{B}$, and $\partial_x^\ell\vc{E}$ and adding the resulting equations, we get
\begin{equation}
\partial_t{\mathcal E}^\ell+\Div \vec {\mathcal F}^\ell
+\nu\left(\partial_x^\ell\vu\right)^2
={\mathcal R}^\ell+{\mathcal S}^\ell,
\label{EFRS}
\end{equation}
where
\[
{\mathcal E}^\ell:=
\frac{1}{2}\ \left(\partial_x^\ell\vu\right)^2
+\frac{1}{2}\ \frac{p_\vr}{\vr}\left(\partial_x^\ell\vr\right)^2
+\frac{1}{2}\ \frac{C_v}{\vt}\left(\partial_x^\ell\vt\right)^2
+\frac{1}{2}\ \frac{1}{4a\vr E_r}\left(\partial_x^\ell E_r\right)^2
+\frac{1}{2}\ \left(\partial_x^\ell\vec E\right)^2
+\frac{1}{2}\ \left(\partial_x^\ell\vec B\right)^2
\]
\[
\vec{\mathcal F}^\ell:=
\left(
\frac{p_\vr}{\vr}\partial_x^\ell\vr
+\frac{p_\vt}{\vr}\partial_x^\ell\vt
+\frac{1}{3a\vr}\partial_x^\ell E_r
\right)
\partial_x^\ell\vu
+
\frac{1}{2}\ \left(
\left(\partial_x^\ell\vu\right)^2
+ \frac{p_\vr}{\vr}\left(\partial_x^\ell\vr\right)^2
+ \frac{C_v}{\vt}\left(\partial_x^\ell\vt\right)^2
+ \frac{1}{4a\vr E_r}\left(\partial_x^\ell E_r\right)^2
\right)\vu
\]
\[
{\mathcal R}^\ell:=
\frac{1}{2} \left[\frac{p_\vr}{\vr^2}\right]_t \left(\partial_x^\ell\vr\right)^2
+\frac{1}{2} \left[\frac{C_v}{\vt}\right]_t \left(\partial_x^\ell\vt\right)^2
+\frac{1}{2} \left[\frac{1}{4a\vr E_r}\right]_t \left(\partial_x^\ell E_r\right)^2
\]
\[
+\frac{1}{2} \Div\left(\frac{p_\vr}{\vr^2}\vu\right)\left(\partial_x^\ell\vr\right)^2
+\frac{1}{2} \Div\vu\left(\partial_x^\ell\vu\right)^2
+\frac{1}{2} \Div\left(\frac{C_v}{\vt}\vu\right)\left(\partial_x^\ell\vt\right)^2
+\frac{1}{2} \Div\left(\frac{1}{4a\vr E_r}\vu\right)\left(\partial_x^\ell E_r\right)^2
\]
\[
+\Grad\left(\frac{p_\vr}{\vr}\right)\partial_x^\ell\vr\ \partial_x^\ell\vu
+\Grad\left(\frac{p_\vt}{\vr}\right)\partial_x^\ell\vt\ \partial_x^\ell\vu
+\Grad\left(\frac{1}{3a\vr}\right)\partial_x^\ell E_r\ \partial_x^\ell\vu
\]
\[
+\frac{p_\vr}{\vr^2}\partial_x^\ell\vr\ F_1^\ell
+\partial_x^\ell\vu\ F_2^\ell
+\frac{C_v}{\vt}\partial_x^\ell\vt\ F_3^\ell
+\partial_x^\ell E_r\ F_4^\ell
+\overline\vr \partial_x^\ell\vec E\cdot\partial_x^\ell\vu,
\]
and
\[
{\mathcal S}^\ell:=
-\partial_x^\ell \vu\cdot\partial_x^\ell\left[\vc u\times\left(\vec B-\overline{\vec B}\right)\right]
 - \frac{C_v}{\vt}\partial_x^\ell \vt\ \partial_x^\ell\left[\frac{\sigma_a}{\vr C_v}\left(a\vt^4-E_r\right)\right]
\]
\[
-\frac{1}{4a\vr E_r}\partial_x^\ell E_r\ \partial_x^\ell\left[\sigma_a\left(E_r-a\vt^4\right)\right]
+\partial_x^\ell\vec E\ \partial_x^\ell[(\vr-\overline\vr)\vec u]
\]
Integrating (\ref{EFRS}) on space, one gets
\[
\partial_t\int_{{\mathbb R}^3}{\mathcal E}^\ell dx
+\left\|\partial_x^\ell\vu\right\|_{L^2}^2
\leq
\int_{{\mathbb R}^3}(|{\mathcal R}^\ell|+|{\mathcal S}^\ell|)\ dx.
\]
Integrating now with respect to $t$ and summing on $\ell$ with $|\ell|\leq d$, we get
\[
\left\|\partial_x V(t) \right\|_{H^{d-1}}^2 + \int_0^t \|\partial_x \vu(\tau)\|^2_{H^{d-1}}
\ d\tau
\leq C_0 \left\|\partial_x V_0 \right\|_{H^{d-1}}^2
+C_0\sum_{|\ell|=1}^d\int_{{\mathbb R}^3}(|{\mathcal R}^\ell|+|{\mathcal S}^\ell|)\ dx.
\]
Observing that
$|\partial_t\vr|\leq C|\partial_x\vr|$, $|\partial_t\vt|\leq C(|\partial_x\vr|+|\partial_x\vt|+|\partial_x E_r||\Delta\vt|)$
 and $|\partial_t E_r|\leq C(|\partial_x\vr|+|\partial_x\vt|+|\partial_x E_r|)$,
and that, using commutator estimates (see Moser-type calculus inequalities in \cite{M})
\[
\|(F_1^\ell,F_2^\ell,F_3^\ell,F_4^\ell)\|_{L^2}
\leq
\|\partial_x(\vr-\overline\vr,\vu,\vt-\overline\vt,E_r-\overline E_r)\|_{L^\infty}
\|\partial_x^\ell(\vr-\overline\vr,\vu,\vt-\overline\vt,E_r-\overline E_r)\|_{L^2}^2,
\]
 we see that
\[
|{\mathcal R}^\ell|
\leq
C\left(\|\partial_x\vr\|_{L^\infty}+\|\partial_x\vu\|_{L^\infty}+\|\partial_x\vt\|_{L^\infty}
+\|\partial_x E_r\|_{L^\infty}
\right)
\|\partial_x^\ell(\vr-\overline\vr,\vu,\vt-\overline\vt,E_r-\overline E_r)\|_{L^2}^2.
\]
Then integrating with respect to time
\[
\int_0^t
|{\mathcal R}^\ell(\tau)|\ d\tau
\]
\[
\displaystyle{
\leq
C\sup_{0\leq\tau\leq t}
\left\{
\|\partial_x\vr\|_{L^\infty}+\|\partial_x\vu\|_{L^\infty}+\|\partial_x\vt\|_{L^\infty}
+\|\partial_x E_r\|_{L^\infty}
\right\}
\int_0^t \|\partial_x^\ell(\vr-\overline\vr,\vu,\vt-\overline\vt,E_r-\overline E_r)\|_{L^2}^2 d\tau
}
\]
\[
\leq CE(t)D^2(t),
\]
for any $|\ell|\leq d$. In the same stroke, we estimate
\[
\displaystyle{
|{\mathcal S}^\ell|
\leq
 C\|\partial_x^\ell \vu\|_{L^2}^2 \left\| \partial_x^\ell\left[\vc u\times\left(\vec B-\overline{\vec B}\right)\right]\right\|_{L^2}^2
+C\|\partial_x^\ell \vt\|_{L^2}^2 \left\| \partial_x^\ell\left[\frac{\sigma_a}{\vr C_v}\left(a\vt^4-E_r\right)\right]\right\|_{L^2}^2
}
\]
\[
\displaystyle{
+C\|\partial_x^\ell E_r\|_{L^2}^2 \left\| \partial_x^\ell\left[\sigma_a\left(E_r-a\vt^4\right)\right]\right\|_{L^2}^2
+C\|\partial_x^\ell\vec E \|_{L^2}^2 \left\| \partial_x^\ell[(\vr-\overline\vr)\vec u]\right\|_{L^2}^2.
}
\]
Then we get
\[
\displaystyle{
|{\mathcal S}^\ell|
\leq
C\|\vec B-\overline{\vec B}\|_{L^\infty}
\|\partial_x^\ell\vu\|_{L^2}^2
}
\]
\[
\displaystyle{
+
C\|\left(\vr-\overline\vr,\vu,\vt-\overline\vt,E_r-\overline E_r\right)\|_{L^\infty}
\|\partial_x^\ell\left(\vr-\overline\vr,\vu,\vt-\overline\vt,E_r-\overline E_r\right)\|_{L^2}
\|\partial_x^\ell(\vec B,\vec E)\|_{L^\infty}
}
\]
\[
\displaystyle{
+
C\left(\|\partial_x\vr\|_{L^\infty}
+\|\partial_x\vu\|_{L^\infty}
+\|\partial_x\vt\|_{L^\infty}
+\|\partial_x E_r\|_{L^\infty}
\right)
\|\partial_x^\ell(\vr-\overline\vr,\vu,\vt-\overline\vt,E_r-\overline E_r)\|_{L^2}^2.
}
\]
Then integrating with respect to time
\[
\displaystyle{
\int_0^t
|{\mathcal S}^\ell(\tau)|\ d\tau
\leq
C\sup_{0\leq\tau\leq t}\|(\vec B-\overline{\vec B})(\tau)\|_{L^\infty}\int_0^t \|\partial_x^\ell\vu(\tau)\|_{L^2}^2 d\tau
}
\]
\[
\displaystyle{
+C\sup_{0\leq\tau\leq t}\|\partial_x^\ell(\vec B,\vec E)(\tau)\|_{L^2}
\int_0^t \|\left(\vr-\overline\vr,\vu,\vt-\overline\vt,E_r-\overline E_r\right)(\tau)\|_{L^\infty}
\|\partial_x^\ell\left(\vr-\overline\vr,\vu,\vt-\overline\vt,E_r-\overline E_r\right)(\tau)\|_{L^2} d\tau
}
\]
\[
\displaystyle{
+C\sup_{0\leq\tau\leq t}\left\{\|\partial_x\vr\|_{L^\infty}+\|\partial_x\vu\|_{L^\infty}+\|\partial_x\vt\|_{L^\infty}
+\|\partial_x E_r\|_{L^\infty}
(\tau)\right)
}
\]
\[
\times\int_0^t \|\partial_x^\ell(\vr-\overline\vr,\vu,\vt-\overline\vt,E_r-\overline E_r)\|_{L^2}^2 d\tau
\]
\[
\leq C\left(E(t)D^2(t)+F(t)I(t)D(t)\right),
\]
for any $|\ell|\leq d$.
\hfill$\square$

\medskip

The above results, together with (\ref{eq:linftyhd_estimat}), allow to derive the following energy bound:

\begin{Corollary}
\label{pr:fin_liftyhd_estimat}
  Assume that the assumptions of Proposition~\ref{pr:linftyh_estimat}
  are satisfied. Then
\begin{equation}
\label{eq:fin_lifty_hd}
 \left\|(V - \overline V)(t) \right\|_{H^d}^2
 + \int_0^t
 \| \vu(\tau)\|^2_{H^d}\ d\tau
\leq
 C\left\|(V - \overline V)(0) \right\|_{H^d}^2
+C\left(E(t)D(t)^2+F(t)I(t)D(t)\right).
  \end{equation}
\end{Corollary}
\vskip0.25cm
Our goal is now to derive bounds for the integrals in the right-hand and left-hand sides of equation \eqref{eq:fin_lifty_hd}.
For that purpose we adapt the results of Ueda, Wang and Kawashima \cite{UWK}.

\newtheorem{luwk1}[l1]{Lemma}
\begin{luwk1}
\label{luwk1}
Under the same assumptions as in Theorem~\ref{th:existence_Maxwell}, and supposing that $d\geq 3$,
we have the following estimate for any $\varepsilon>0$
\[
\int_0^t\left(\|\left(\vr-\overline\vr,\vt-\overline\vt,E_r-\overline E_r\right)(\tau)\|^2_{H^{d}}+ \left\| \vec E(\tau) \right\|^2_{H^{d-1}}\right)\ d\tau
\]
\begin{equation}
  \label{eq:estimate_1}
 \leq
\varepsilon \int_0^t \left\|\partial_x \vec B(\tau) \right\|^2_{H^{d-2}}\ d\tau
+C_{\varepsilon} \left\{ \|V_0-\overline{V}\|^2_{H^{d-1}}+E(t)D(t)^2+F(t)I(t)D(t)\right\}.
\end{equation}
\end{luwk1}
\noindent {\bf Proof:} We linearize the principal part of the system
\eqref{j1bis}-\eqref{j2bis}-\eqref{j3bis} as follows

\bFormula{k1bis}
\partial_t \vr + \overline\vr\ \Div\vu=g_1,
\eF
\bFormula{k2bis}
\partial_t \vu+\overline a_1\ \Grad \vr+\overline a_2\ \Grad \vt + \overline a_3\ \Grad E_r +\vec E+\vc u\times\overline{\vec B}+\nu\vu = g_2,
\eF
\bFormula{k3bis}
\partial_t \vt+\overline b_1\ \Div \vu+\overline b_2(\vt-\overline\vt)
=g_3,
\eF
\bFormula{k4bis}
\partial_t E_r + \overline c_1\ \Div \vu +\overline c_3(E_r-\overline E_r)
 =g_4,
\eF
\bFormula{k5bis}
\partial_t \vc{B} + \Curl \vc{E} = 0,\
\eF
\bFormula{k6bis}
\partial_t \vc{E} - \Curl \vc{B} -\overline \vr\ \vu =g_5,
\eF
with coefficients
\[
a_1(\vr,\vt)=\frac{p_\vr}{\vr},\ \
a_2(\vr,\vt)=\frac{p_\vt}{\vr},\ \
a_3(\vr,\vt)=\frac{1}{3\vr},\ \
\overline a_j=a_j(\overline\vr,\overline\vt),
\]
\[
b_1(\vr,\vt)=\frac{\vt p_\vt}{\vr C_v},\ \
b_2(\vr,\vt,E_r)=\frac{a\sigma_a}{\vr C_v}(\vt^2+\overline\vt^2)(\vt+\overline\vt),\ \
b_3(\vr,\vt,E_r)=\frac{a\sigma_a}{\vr C_v},\ \
\overline b_j= b_j(\overline\vr,\overline\vt),
\]
\[
c_1(\vr,\vt,E_r)=\frac{4}{3}E_r,\ \
c_2(\vr,\vt,E_r)=a\sigma_a(\vt^2+\overline\vt^2)(\vt+\overline\vt),\ \
c_3(\vr,\vt,E_r)=\sigma_a,\ \
\overline c_j= c_j(\overline\vr,\overline\vt),
\]
and sources
\[
g_1:=-\{\vu\cdot\Grad\vr+(\vr-\overline\vr)\Div\vu\},
\]
\[
g_2:=-\left\{(\vu\cdot\Grad)\vu+(a_1-\overline a_1)\Grad \vr+(a_2-\overline a_2)\Grad \vt + (a_3-\overline a_3)\Grad E_r +\vc u\times(\vec B-\overline{\vec B})\right\},
\]
\[
g_3:=-\left\{(\vu\cdot\Grad)\vt+ (b_1-\overline b_1)\Div \vu+(b_2-\overline b_2)(\vt-\overline\vt)
+b_3(E_r-\overline E_r)
\right\},
\]
\[
g_4:=-\left\{(\vu\cdot\Grad)E_r+(\overline c_1-c_1)\Div \vu
+c_2(\vt-\overline\vt)
+(c_3-\overline c_3)(E_r-\overline E_r)  \right\},
\]
and
\[
g_5=(\vr-\overline \vr) \vu.
\]
Multiplying
 (\ref{k1bis}) by $-\overline a_1 \Div\vu$,
(\ref{k2bis}) by $\overline a_1\ \Grad \vr+\overline a_2\ \Grad \vt + \overline a_3\ \Grad E_r +\vec E$,
(\ref{k3bis}) by $-\overline{a}_2\Div \vu
+\vt-\overline\vt$,
(\ref{k4bis}) by $-\overline{a}_3\Div \vu
+E_r-\overline E_r$,
(\ref{k5bis}) by $1$,
(\ref{k6bis}) by $\vu$ and summing up, we get
\[
\overline a_1( \Grad \vr\ \vu_t-\vr_t\ \Div \vu)
+
\overline a_2( \Grad \vt\ \vu_t-\vt_t\ \Div \vu)
+
\overline a_3( \Grad E_r\ \vu_t-(E_r)_t\ \Div \vu)
+\vec E\vu_t+\vec E_t\vu
\]
\[
+\left\{\frac{1}{2}\left[(\vt-\overline\vt)^2+(E_r-\overline E_r)^2\right]\right\}_t
\]
\[
+(\overline a_1 \Grad \vr+ \overline a_2 \Grad \vt+\overline a_3 \Grad E_r+\vec E)^2
+(\overline a_1\ \Grad \vr+\overline a_2\ \Grad \vt + \overline a_3\ \Grad E_r +\vec E)(\vu\times\overline{\vec B}+\nu\vu)
\]
\[
+\overline b_2(\vt-\overline\vt)^2+\overline c_3(E_r-\overline E_r)^2
\]
\[
+\overline b_1 (\vt-\overline\vt)  \Div \vu+\overline c_1 (E_r-\overline E_r)  \Div \vu
\]
\[
+(\overline a_3\overline c_2-\overline a_2\overline b_2)(\vt-\overline\vt)\Div\vu
+(\overline a_2\overline b_3-\overline a_3\overline c_3)(E_r-\overline E_r)\Div\vu
\]
\begin{equation}
-\vu\ \Curl\vec B-\overline\vr\vu^2
-(\Div\vu)^2\left[ \overline a_1+\overline a_2+\overline a_3\right]
=G_1^0,
\label{ident1}
\end{equation}
where
\[
G_1^0:=
-\overline a_1 g_1  \Div \vu
+[\overline a_1 \Grad \vr+ \overline a_2 \Grad \vt+\overline a_3 \Grad E_r+\vec E]g_2
-[\overline a_2 +\vt-\overline\vt]  \Div \vu g_3
-[\overline a_3 +E_r-\overline E_r]  \Div \vu g_4
+g_5\vu.
\]
Rearranging the left hand side of (\ref{ident1}) we get
\begin{equation}
\{H_1^0\}_t+\Div \vec F_1^0+D_1^0
=M_1^0+G_1^0,
\label{ident2}
\end{equation}
where
\[
H_1^0
=-\left[
\overline a_1(\vr-\overline \vr)+\overline a_2(\vt-\overline \vt)+\overline a_3(E_r-\overline E_r)\right]\Div\vu
+\vec E\cdot\vu
+\frac{1}{2}\left[(\vt-\overline\vt)^2+(E_r-\overline E_r)^2\right],
\]
\[
\vec F_1^0
=
\left[
\overline a_1(\vr-\overline \vr)+\overline a_2(\vt-\overline \vt)+\overline a_3(E_r-\overline E_r)\right]\vu_t
-2\left[\overline a_1(\vr-\overline \vr)+\overline a_2(\vt-\overline \vt)+\overline a_3(E_r-\overline E_r)\right]\vec E
\]
\[
+(\overline a_3\overline c_2-\overline a_2\overline b_2+\overline b_1)(\vt-\overline\vt)\vu
+(\overline a_2\overline b_3-\overline a_3\overline c_3+\overline c_1)(E_r-\overline E_r)\vu,
\]
\[
D_1^0
=
\overline a_1^2 |\Grad \vr|^2+ \overline a_2^2 |\Grad \vt|^2+\overline a_3^2 |\Grad E_r|^2+|\vec E|^2
+2\overline a_1(\vr-\overline \vr)^2+\overline b_2(\vt-\overline\vt)^2+\overline c_3(E_r-\overline E_r)^2,
\]
and
\[
M_1^0
=-\left\{
2\overline a_1\overline a_2\Grad \vr\cdot \Grad \vt
+2\overline a_1\overline a_3\Grad \vr\cdot \Grad E_r
+2\overline a_2\overline a_3\Grad \vt\cdot \Grad E_r
\right.
\]
\[
+2\overline a_2(\vr-\overline \vr)(\vt-\overline \vt)
+2\overline a_2(\vr-\overline \vr)(E_r-\overline E_r)
\]
\[
+(\overline a_1\ \Grad \vr+\overline a_2\ \Grad \vt + \overline a_3\ \Grad E_r +\vec E)(\vu\times\overline{\vec B}+\nu\vu)
-\vu\ \Curl\vec B-\overline\vr\vu^2
-(\Div\vu)^2\left[ \overline a_1+\overline a_2+\overline a_3\right]
\]
\[
\left.
-(\overline a_3\overline c_2-\overline a_2\overline b_2+\overline b_1)\Grad\vt\cdot\vu
-(\overline a_2\overline b_3-\overline a_3\overline c_3+\overline c_1)\Grad E_r\cdot\vu
\right\}.
\]
Integrating (\ref{ident2}) over space and using Young's inequality, we find
\[
\frac{d}{dt}\int_{{\mathbb R}^3}H_1^0\ dx
+C\left(
\|
 \vr\|^2_{L^2}
+ \|\Grad \vt\|^2_{L^2}
+\|\Grad E_r\|^2_{L^2}
+\|\vec E\|^2_{L^2}
+\|\vr-\overline \vr\|^2_{L^2}
\right)
\]
\begin{equation}
\leq
\varepsilon \|\partial_x\vec B\|^2_{L^2}
+C_\varepsilon\left(\|\vu\|^2_{H^1}
+\|\vt-\overline \vt\|^2_{H^1}+\|E_r-\overline E_r\|^2_{H^1}\right)
+ \int_{{\mathbb R}^3}|G_1^0|\ dx.
\label{ident3}
\end{equation}
In fact one obtains in the same way estimates for the derivatives of $ V$.

Namely, applying $\partial_x^\ell$ to the system (\ref{k1bis}-\ref{k6bis}), we get

\begin{equation}
\{H_1^\ell\}_t+\Div \vec F_1^\ell+D_1^\ell
=M_1^\ell+G_1^\ell,
\label{identl}
\end{equation}
where
\[
H_1^\ell
=-\left[
\overline a_1\partial_x^\ell (\vr-\overline\vr)+\overline a_2\partial_x^\ell(\vt-\overline\vt)+\overline a_3\partial_x^\ell( E_r-\overline E_r)\right]\Div\partial_x^\ell\vu
+\partial_x^\ell\vec E\cdot\partial_x^\ell\vu
\]
\[
+\frac{1}{2}[(\partial_x^\ell\vt)^2+(\partial_x^\ell E_r)^2],
\]
\[
\vec F_1^\ell
=
\left[
\overline a_1\partial_x^\ell (\vr-\overline \vr)+\overline a_2\partial_x^\ell (\vt-\overline \vt)+\overline a_3\partial_x^\ell (E_r-\overline E_r)\right]\vu_t
\]
\[
+(\overline a_3\overline c_2-\overline a_2\overline b_2+\overline b_1)\partial_x^\ell\vt\ \partial_x^\ell\vu
+(\overline a_2\overline b_3-\overline a_3\overline c_3+\overline c_1)\partial_x^\ell E_r\ \partial_x^\ell\vu,
\]
\[
-2\left[\overline a_1\partial_x^\ell (\vr-\overline \vr)+\overline a_2\partial_x^\ell (\vt-\overline \vt)+\overline a_3\partial_x^\ell (E_r-\overline E_r)\right]\partial_x^\ell \vec E
+\partial_x^\ell \vu \times \partial_x^\ell (\vec B-\overline{\vec B}) ,
\]
\[
D_1^\ell
=
\overline a_1^2 |\Grad \partial_x^\ell \vr|^2+ \overline a_2^2 |\partial_x^\ell \Grad \vt|^2+\overline a_3^2 |\partial_x^\ell \Grad E_r|^2+|\partial_x^\ell \vec E|^2
+2\overline a_1\left(\partial_x^\ell (\vr-\overline \vr)\right)^2+\overline b_2(\partial_x^\ell\vt)^2+\overline c_3(\partial_x^\ell E_r)^2,
\]
\[
M_1^\ell
=-\left\{
2\overline a_1\overline a_2\Grad \partial_x^\ell \vr\cdot \Grad \partial_x^\ell \vt
+2\overline a_1\overline a_3\Grad \partial_x^\ell \vr\cdot \Grad \partial_x^\ell E_r
+2\overline a_2\overline a_3\Grad \partial_x^\ell \vt\cdot \Grad \partial_x^\ell E_r
\right.
\]
\[
+2\overline a_2\partial_x^\ell (\vr-\overline \vr)\partial_x^\ell (\vt-\overline \vt)
+2\overline a_2\partial_x^\ell (\vr-\overline \vr)\partial_x^\ell (E_r-\overline E_r)
+(\overline a_1\ \Grad \partial_x^\ell \vr+\overline a_2\ \Grad \partial_x^\ell \vt + \overline a_3\ \Grad \partial_x^\ell E_r +\partial_x^\ell \vec E)
(\partial_x^\ell \vu\times\overline{\vec B}+\nu\partial_x^\ell \vu)
\]
\[
-(\overline a_3\overline c_2-\overline a_2\overline b_2+\overline b_1)\Grad\partial_x^\ell\vt\cdot \partial_x^\ell\vu
-(\overline a_2\overline b_3-\overline a_3\overline c_3+\overline c_1)\Grad \partial_x^\ell E_r\cdot \partial_x^\ell\vu
\]
\[
\left.
-\Curl\partial_x^\ell \vu\ \partial_x^\ell (\vec B-\overline{\vec B})-\overline\vr\left(\partial_x^\ell \vu\right)^2
-(\Div\partial_x^\ell \vu)^2\left[ \overline a_1+\overline a_2+\overline a_3\right]
\right\},
\]
and
\[
G_1^\ell
=
-\overline a_1 \partial_x^\ell g_1  \Div \partial_x^\ell \vu
+[\overline a_1 \Grad \partial_x^\ell \vr+ \overline a_2 \Grad \partial_x^\ell \vt+\overline a_3 \Grad \partial_x^\ell E_r+\partial_x^\ell \vec E]\partial_x^\ell g_2
\]
\[
-\overline a_2 \partial_x^\ell g_3  \Div \partial_x^\ell \vu
-\overline a_3 \partial_x^\ell g_4  \Div \partial_x^\ell \vu
+\partial_x^\ell g_5\partial_x^\ell \vu
+\partial_x^\ell g_3 \partial_x^\ell \vt
+\partial_x^\ell g_4 \partial_x^\ell E_r.
\]
Integrating (\ref{identl}) over space and time, we find
\[
\int_{{\mathbb R}^3}H_1^\ell(t)\ dx
-\int_{{\mathbb R}^3}H_1^\ell(0)\ dx
\]
\[
+C\int_0^t
\left(
\|\Grad \partial_x^\ell \vr\|^2_{L^2}
+ \|\Grad \partial_x^\ell \vt\|^2_{L^2}
+\|\Grad \partial_x^\ell E_r\|^2_{L^2}
+\|\partial_x^\ell \vec E\|^2_{L^2}
\right)\ d\tau
\]
\[
+C\int_0^t
\left(
+\|\partial_x^\ell (\vr-\overline \vr)\|^2_{L^2}
+\|\partial_x^\ell (\vt-\overline \vt)\|^2_{L^2}
+\|\partial_x^\ell (E_r-\overline E_r)\|^2_{L^2}
\right)\ d\tau
\]
\[
\leq
\varepsilon\int_0^t  \|\partial_x^\ell (\vec B-\overline{\vec B})\|^2_{L^2}\ d\tau
+C_\varepsilon\int_0^t
\left(\|\partial_x^\ell \vu\|^2_{H^1}
+\|\partial_x^\ell (\vt-\overline \vt)\|^2_{H^1}+\|\partial_x^\ell (E_r-\overline E_r)\|^2_{H^1}\right)\ d\tau
\]
\begin{equation}
+\int_0^t \int_{{\mathbb R}^3}|G_1^\ell|\ dx\ d\tau.
\label{identlbis}
\end{equation}
Observing that
\[
\left|
\int_{{\mathbb R}^3}H_1^\ell(t)\ dx
\right|
\leq
C\left(
\|\partial_x^\ell (\vr-\overline\vr)\|^2_{L^2}
+\|\partial_x^\ell(\vt-\overline\vt)\|^2_{L^2}
+\|\partial_x^\ell( E_r-\overline E_r)\|^2_{L^2}
+\|\partial_x^\ell\vu\|^2_{H^1}\right),
\]
and summing (\ref{identlbis}) on $\ell$ for $1\leq\ell\leq d-1$, we get
\[
\int_0^t
\left(
\|\left(\vr-\overline\vr,\vt-\overline\vt, E_r-\overline E_r\right)(\tau)\|^2_{H^d}+\| \vec E(\tau)\|^2_{H^{d-1}}\right)\ d\tau
\leq
C_\varepsilon \|(V-\overline V)(0)\|^2_{H^d}
\]
\begin{equation}
+\varepsilon \int_0^t  \|\partial_x\vec B(\tau)\|^2_{H^{d-2}}\ d\tau
+C_\varepsilon \left(E(t)D^2(t)+F(t)I(t)D(t)\right)
+\sum_{|\ell|=1}^{d-1}\int_0^t \int_{{\mathbb R}^3}|G_1^\ell(\tau)|\ dx\ d\tau,
\label{identlter}
\end{equation}
where we used Corollary \ref{pr:fin_liftyhd_estimat}.

Let us estimate the last integral in (\ref{identlbis}). we have
\begin{equation}
\left\{
\begin{array}{c}
\|\partial_x^\ell g_1\|_{L^2}
\leq
C\|(\vr-\overline\vr,\vu,\vt-\overline\vt, E_r-\overline E_r)\|_{L^\infty} \|\partial_x^{\ell+1}(\vr,\vu)\|_{L^2},\\\\
\|\partial_x^\ell g_2\|_{L^2}
\leq
C\|(\vr-\overline\vr,\vu,\vt-\overline\vt, E_r-\overline E_r)\|_{L^\infty} \|\partial_x^{\ell+1}(\vr,\vu)\|_{L^2}\\\\
+C\|\vec B-\overline{\vec B}\|_{L^\infty} \|\partial_x^\ell\vu\|_{L^2}
+C\|\partial_x^\ell(\vec B-\overline{\vec B})\|_{L^2}\|\vu\|_{L^\infty},\\\\
\|\partial_x^\ell g_3\|_{L^2}
\leq
C\|(\vr-\overline\vr,\vu,\vt-\overline\vt, E_r-\overline E_r)\|_{L^\infty} \|\partial_x^{\ell+1}(\vr,\vu)\|_{L^2}\\\\
+C\|(\vr-\overline\vr,\vu,\vt-\overline\vt, E_r-\overline E_r)\|_{L^\infty} \|\partial_x^{\ell+2}(\vt,E_r)\|_{L^2},\\\\
\|\partial_x^\ell g_4\|_{L^2}
\leq
C\|(\vr-\overline\vr,\vu,\vt-\overline\vt, E_r-\overline E_r)\|_{L^\infty} \|\partial_x^{\ell+1}(\vr,\vu)\|_{L^2}\\\\
+C\|(\vr-\overline\vr,\vu,\vt-\overline\vt, E_r-\overline E_r)\|_{L^\infty} \|\partial_x^{\ell+2}(\vt,E_r)\|_{L^2},\\\\
\|\partial_x^\ell g_5\|_{L^2}
\leq
C\|(\vr-\overline\vr,\vu,\vt-\overline\vt, E_r-\overline E_r)\|_{L^\infty} \|\partial_x^{\ell}(\vr,\vu)\|_{L^2},
\end{array}
\right.
\label{Ineg}
\end{equation}
for $1\leq |\ell|\leq d-1$. Then
\[
\int_0^t \int_{{\mathbb R}^3}|G_1^\ell(\tau)|\ dx\ d\tau
\leq
C\|\partial_x^{\ell+1}\vu\|_{L^2}\|\partial_x^\ell g_1\|_{L^2}
\]
\[
+C \left( \|\partial_x^{\ell+1}\vr\|_{L^2}+\|\partial_x^{\ell+1}\vt\|_{L^2}+\|\partial_x^{\ell+1}E_r\|_{L^2}+\|\partial_x^{\ell}\vec E\|_{L^2}\right)
\|\partial_x^\ell g_2\|_{L^2}
\]
\[
+C\|\partial_x^{\ell+1}\vu\|_{L^2}\|\partial_x^\ell g_3\|_{L^2}
+C\|\partial_x^{\ell+1}\vu\|_{L^2}\|\partial_x^\ell g_4\|_{L^2}
+C\|\partial_x^{\ell}\vu\|_{L^2}\|\partial_x^\ell g_5\|_{L^2}.
\]
Plugging bounds (\ref{Ineg}) into this last inequality gives
\[
\sum_{|\ell|=1}^{d-1}\int_0^t \int_{{\mathbb R}^3}|G_1^\ell(\tau)|\ dx\ d\tau
\leq
CE(t)D^2(t),
\]
which ends the proof of Lemma \ref{luwk1}.  \hfill $\square$
\vskip0.5cm

Finally we check from \cite{UWK} (see Lemma 4.4) that the following result for the Maxwell's system holds true for our system with a similar proof

\newtheorem{luwk2}[l1]{Lemma}
\begin{luwk2}
\label{luwk2}
Under the same assumptions as in Theorem~\ref{th:existence_Maxwell}, and supposing that $d\geq 3$,
 for any $\varepsilon>0$ the following estimate (here, we set $V= (\vr,\vu,\vt,E_r,\vec B,\vec E)^T$) holds
\begin{equation}
  \label{eq:estimate_2}
 \int_0^t \left\|\partial_x \vec B(\tau) \right\|^2_{H^{s-2}}\ d\tau
 \leq
C\|V_0-\overline{V}\|^2_{H^{s-1}}
+C\int_0^t \left\|\partial_x \vec E(\tau) \right\|^2_{H^{s-2}}\ d\tau
+C(E(t)D(t)^2+F(t)I(t)D(t)).
\end{equation}
\end{luwk2}
\noindent {\bf Proof:} Applying $\partial_x^\ell$ to (\ref{j5bis}) and (\ref{j6bis}), multiplying respectively by
 $-\Curl\partial_x^\ell\vec B$, (\ref{j6bis}) by $\Curl\partial_x^\ell\vec E$ and adding the resulting equations, we get
\[
-\left(\partial_x^\ell\vec E\cdot\Curl\partial_x^\ell\vec B\right)_t+|\Curl\partial_x^\ell\vec B|^2
-\Div\left(\partial_x^\ell\vec E\times\partial_x^\ell\vec B_t\right)=M_2^\ell+G_2^\ell,
\]
where
\[
M_2^\ell=-\overline\vr\partial_x^\ell\vu\cdot\Curl\partial_x^\ell\vec B+|\Curl\partial_x^\ell\vec E|^2,
\]
and
\[
G_2^\ell=-\partial_x^\ell\left((\vr-\overline\vr)\vu\right)\cdot\Curl\partial_x^\ell\vec B.
\]
Integrating in space we get
\[
-\frac{d}{dt}\int_{{\mathbb R}^3}\partial_x^\ell\vec E\cdot\Curl\partial_x^\ell\vec B\ dx
+C\|\Curl\partial_x^\ell\vec B\|^2_{L^2}
\leq
\|\Curl\partial_x^\ell\vec E\|^2_{L^2}
+\|\partial_x^\ell\vu\|^2_{L^2}
+\int_{{\mathbb R}^3}|G_2^\ell|\ dx.
\]
Integrating on time and summing for $1\leq |\ell|\leq d-2$, we have
\[
\int_0^t \|\partial_x\vec B\|^2_{H^{d-2}} dt
\leq
C\|(V-\overline V)(t)\|_{H^{d-1}}+C\|(V-\overline V)(0)\|_{H^{d-1}}
+C\int_0^t \|\partial_x\vec E\|^2_{H^{d-2}} dt
\]
\[
+C\int_0^t \|\vu\|^2_{H^{d-2}} dt
+C\sum_{|\ell|=0}^{d-2}\int_{{\mathbb R}^3}|G_2^\ell(\tau)|\ dx\ d\tau
\]
\[
\leq
C\|(V-\overline V)(0)\|_{H^{d-1}}
+C\int_0^t \|\partial_x\vec E|^2_{H^{d-2}} dt
+C(E(t)D(t)^2+F(t)I(t)D(t)),
\]
where we used the bound
${\displaystyle \sum_{|\ell|=1}^{d-1}\int_0^t \int_{{\mathbb R}^3}|G_2^\ell(\tau)|\ dx\ d\tau
\leq
CE(t)D^2(t)}$, obtained in the same way as in the proof of Lemma \ref{luwk1}, which ends the proof of Lemma \ref{luwk2}.  \hfill$\square$
\medskip

We are now in position to conclude with the proofs of Theorems \ref{th:existence_Maxwell} and \ref{th:asymptotics_Maxwell}.

\subsection{Proof of Theorem~\ref{th:existence_Maxwell}:}
\vskip0.25cm
 We first point out that local existence for the hyperbolic system \eqref{j1bis}-\eqref{j6bis} may be proved using standard fixed-point methods.
 We refer to \cite{M} for the proof.
\vskip0.25cm
Now plugging (\ref{eq:estimate_2}) into (\ref{eq:estimate_1}) with $\varepsilon$ small enough, we get
\[
\int_0^t\left(\|\left(\vr-\overline\vr,\vt-\overline\vt,E_r-\overline E_r\right)\|^2_{H^{d}}+ \left\| \vec E(\tau) \right\|^2_{H^{d-1}}\right)\ d\tau
\]
\begin{equation}
\label{E1}
 \leq
C \left\{ \|V_0-\overline{V}\|^2_{H^{d-1}}+E(t)D(t)^2+F(t)I(t)D(t)\right\}.
\end{equation}
Putting this last estimate into (\ref{eq:estimate_2}) we find
\begin{equation}
  \label{E2}
 \int_0^t \left\|\partial_x \vec B(\tau) \right\|^2_{H^{s-2}}\ d\tau
 \leq
C\|V_0-\overline{V}\|^2_{H^{s}}
+C(E(t)D(t)^2+F(t)I(t)D(t)).
\end{equation}
Then from (\ref{eq:estimate_1}), (\ref{E1}) and (\ref{E2}) we get
\[
 \|(V-\overline{V})(t)\|^2_{H^{d}}
+\int_0^t\left(\|\left(\vr-\overline\vr,\vt-\overline\vt,E_r-\overline E_r\right)(\tau)\|^2_{H^{d}}
+ \left\| \vec E(\tau)\right\|^2_{H^{d-1}}+\left\|\partial_x \vec B(\tau) \right\|^2_{H^{d-2}}
\right)\ d\tau
\]
\[
\leq
C\|V_0-\overline{V}\|^2_{H^{d}}
+C(E(t)D(t)^2+F(t)I(t)D(t)),
\]
or equivalently
\[
F(t)^2+D(t)^2\leq
C\|V_0-\overline{V}\|^2_{H^{d}}
+C(E(t)D(t)^2+F(t)I(t)D(t)).
\]
Now observing that, provided that $d\geq 3$ one has $\|(V-\overline V)(t)\|_{H^d}\leq E(t)\leq C F(t)$, and that,
 provided that $d\geq 2$ one has $I(t)\leq C D(t)$, for some positive constant $C$, we see that
\[
F(t)^2+D(t)^2\leq
C\|V_0-\overline{V}\|^2_{H^{d}}
+CF(t)D(t)^2.
\]
In order to prove global existence, we argue by contradiction, and assume that $T_c>0$ is the maximum time existence. Then, we necessarily
have
\begin{displaymath}
  \lim_{t\to T_c} N(t) = +\infty,
\end{displaymath}
where $N(t)$ is defined by
\[
N(t):=(F(t)^2+D(t)^2)^{1/2}.
\]
 We are thus reduced to prove that $N$ is bounded. For this purpose, we use the argument used
in \cite{BDN}. After the previous calculation, we have
\begin{equation}\label{eq:12bis}
\forall T\in [0,T_c], \quad  N(t)^2 \leq C\left(\|V_0-\overline{V}\|^2_{H^{d}} + N(t)^3\right).
\end{equation}
Hence, setting $\|V_0-\overline{V}\|_{H^{d}} = \varepsilon$, we have
\begin{equation}\label{eq:11bis}
  \frac{N(t)^2}{\varepsilon^2 + N(t)^3} \leq C.
\end{equation}
Studying the variation of $\phi(N) = N^2/\left(\varepsilon^2 +
  N^3\right)$, we see that $\phi'(0) = 0$, that $\phi$ is
increasing on the interval $\left[0, \left(2\varepsilon^2\right)^{1/3}\right]$ and decreasing on the
interval $\left[\left(2\varepsilon^2\right)^{1/3},+\infty
\right)$. Hence,
\begin{displaymath}
  \max \phi = \phi\left(\left(2\varepsilon^2\right)^{1/3} \right) =
  \frac 13 \left(\frac 2 \varepsilon\right)^{2/3}.
\end{displaymath}
Hence we can
choose $\varepsilon$ small enough to have $\phi(N)\leq C$ for all
$N\in [0,N^*]$, where $N^*>0$, we see that $N\leq N^*$, which contradicts \eqref{eq:12bis}.\hfill $\square$

\section{Large time behaviour}

We have the following analogous of Proposition \ref{pr:linftyh_estimat} for time derivatives

\begin{Corollary}
  \label{tdestim}
Let the assumptions of Theorem~\ref{th:existence_Maxwell} be
satisfied and consider the solution $V:=(\vr,\vu,\vt,E_r,\vec B,\vec E)$ of system
\eqref{j1bis}-\eqref{j2bis}-\eqref{j3bis} on $[0,t]$, for some $t>0$.
Then, one gets for a constant $C_0>0$
\begin{equation}
  \label{tdenergy}
 \left\|\partial_t V(t)\right\|_{H^{d-1}}^2
+ \int_0^t
\left(
\left\|\ \partial_t\left(\rho,\vu,\vt,E_r \right)(\tau)\right\|_{H^{d-1}}^2
+\| \partial_t(\vec B,\vec E)(\tau)\|_{H^{d-2}}^2
\right)\ d\tau
 \leq C_0\left\|V_0 -\overline V \right\|_{H^d}^2.
\end{equation}
\end{Corollary}
\noindent {\bf Proof:} Using System (\ref{sys}) we see that
\[
\left\|\partial_t V\right\|_{H^{d-1}}
\leq
C\left\| V-\overline V\right\|_{H^{d}},
\]
\[
\left\|\ \partial_t\left(\rho,\vu,\vt,E_r \right)\right\|_{H^{d-1}}
\leq
\| \partial_x(\rho,\vu,\vt,E_r,\vec B,\vec E )\|_{H^{d-1}}
+C\|(\rho,\vu,\vt,E_r,\vec B,\vec E )\|_{H^{d-1}},
\]
and
\[
\| \partial_t(\vec B,\vec E)\|_{H^{d-2}}
\leq
\| \partial_x(\vec B,\vec E)\|_{H^{d-2}}
+C\left\|\ \vu\right\|_{H^{d-1}}.
\]
Then for $d\geq 3$, using the uniform estimate $\left\|V -\overline V \right\|_{H^d}^2\leq C$
 of Theorem~\ref{th:existence_Maxwell}, we get estimate (\ref{tdenergy}).  \hfill$\square$

\subsection{Proof of Theorem~\ref{th:asymptotics_Maxwell}:}
\vskip0.25cm
Using Corollary \ref{tdestim}, we get
\[
\int_0^\infty\left|\frac{d}{dt}\|(\vr-\overline\vr,\vu,\vt-\overline\vt,E_r-\overline E_r)(t)\|_{H^{d-1}}\right|\ dt
\]
\[
\leq
2\int_0^\infty\|(\vr-\overline\vr,\vu,\vt-\overline\vt,E_r-\overline E_r)(t)\|_{H^{d-1}}
\|\partial_t(\vr,\vu,\vt,E_r)(t)\|_{H^{d-1}}
\ dt
\]
\[
\leq
\int_0^\infty\|(\vr-\overline\vr,\vu,\vt-\overline\vt,E_r-\overline E_r)(t)\|_{H^{d-1}}^2
+\|\partial_t(\vr,\vu,\vt,E_r)(t)\|_{H^{d-1}}^2\ dt
\leq C_0\left\|V_0 -\overline V \right\|_{H^d}^2.
\]
This implies that
\[
t\to\|(\vr-\overline\vr,\vu,\vt-\overline\vt,E_r-\overline E_r)(t)\|_{H^{d-1}}^2\in L^1(0,\infty)\ \ \mbox{and}\ \
t\to \frac{d}{dt}\|(\vr-\overline\vr,\vu,\vt-\overline\vt,E_r-\overline E_r)(t)\|_{H^{d-1}}\in L^1(0,\infty),
\]
 and then
\[
\|(\vr-\overline\vr,\vu,\vt-\overline\vt,E_r-\overline E_r)(t)\|_{H^{d-1}}\to 0,
\]
 when $t\to\infty$.

Now applying Gagliardo-Nirenberg's inequality, and (\ref{eq:energy_estimate_maxwell}) we get
\[
\|(\vr-\overline\vr,\vu,\vt-\overline\vt,E_r-\overline E_r)(t)\|_{W^{d-2,\infty}}
\leq
\|(\vr-\overline\vr,\vu,\vt-\overline\vt,E_r-\overline E_r)(t)\|_{H^{d-2}}^{1/4}
\|\partial_x^2(\vr,\vu,\vt,E_r)(t)\|_{H^{d-2}}^{3/4}.
\]
So
\[
\|(\vr-\overline\vr,\vu,\vt-\overline\vt,E_r-\overline E_r)(t)\|_{W^{d-2,\infty}}\to 0\ \mbox{when}\ t\to\infty.
\]
\vskip0.25cm
Now in the same stroke
\[
t\to\|\vec E(t)\|_{H^{d-1}}^2\in L^1(0,\infty)\ \ \mbox{and}
\ \ t\to \frac{d}{dt}\|\vec E(t)\|_{H^{d-1}}\in L^1(0,\infty),
\]
 and then
\[
\|\vec E(t)\|_{W^{d-1,\infty}}\to 0\ \mbox{when}\ t\to\infty.
\]
\vskip0.25cm
Finally
\[
t\to\|\partial_x\vec B(t)\|_{H^{d-3}}^2\in L^1(0,\infty)\ \ \mbox{and}
\ \ t\to \frac{d}{dt}\|\partial_x\vec B(t)\|_{H^{d-3}}\in L^1(0,\infty),
\]
 then arguing as before
\[
\|(\vec B-\overline{\vec B})(t)\|_{W^{d-4,\infty}}
\leq
\|(\vec B-\overline{\vec B})(t)\|_{H^{d-4}}^{1/4}
\|\partial_x^2\vec B(t)\|_{H^{d-3}}^{3/4}.
\]
So
\[
\|(\vec B-\overline{\vec B})\|_{W^{d-4,\infty}}\to 0\ \mbox{when}\ t\to\infty,
\]
which ends the proof  \hfill$\square$

\bigskip

{\bf Acknowlegment:}
{\it \v S\' arka Ne\v casov\' a acknowledges the support of the GA\v CR (Czech Science Foundation) project 16-03230S in the framework of RVO: 67985840.
 Bernard Ducomet is partially supported by the ANR project INFAMIE (ANR-15-CE40-0011)}


\vskip1cm
\centerline{Xavier Blanc}
 \centerline{Univ. Paris Diderot, Sorbonne Paris Cit\'e,}
\centerline{Laboratoire Jacques-Louis Lions,}
\centerline{UMR 7598, UPMC, CNRS, F-75205 Paris, France}
 \centerline{E-mail: blanc@ann.jussieu.fr}
\vskip0.5cm
\centerline{Bernard Ducomet}
 \centerline{CEA, DAM, DIF}
\centerline{ F-91297 Arpajon, France}
 \centerline{E-mail: bernard.ducomet@cea.fr}
\vskip0.5cm
\centerline{\v S\' arka Ne\v casov\' a}
\centerline{Institute of Mathematics of the Academy of Sciences of the Czech Republic}
\centerline{\v Zitn\' a 25, 115 67 Praha 1, Czech Republic}
\centerline{E-mail: matus@math.cas.cz}

 \end{document}